\documentclass{article}
\usepackage{amsthm}
\usepackage{amsfonts,amssymb,amsmath}
\usepackage{enumitem}
\usepackage{titlesec}
\usepackage{fouriernc}
\usepackage[T1]{fontenc}
\usepackage{chngcntr}
\counterwithin{figure}{section}

\newlist{gcases}{enumerate}{1}
\setlist[gcases,1]{
  label={{\it Case}~{\it \Alph*}.},
  topsep=0ex,
  leftmargin=0in,
  labelsep=.1in,
  itemindent=.7in,
  itemsep=0ex 
}
\newlist{tenumerate}{enumerate}{1}
\setlist[tenumerate,1]{
  label={(\arabic*)},
  topsep=0ex,
  leftmargin=.3in,
  labelsep=.1in,
  itemindent=0in,
  itemsep=0ex
}

\newlist{titemize}{enumerate}{1}
\setlist[titemize,1]{
  label={$\bullet$},
  topsep=0ex,
  leftmargin=.3in,
  labelsep=.1in,
  itemindent=0in,
  itemsep=0ex
}

\titleformat{\section}
  {\normalfont\bfseries}   
  {}                             
  {0pt}                          
  {Section \thesection\quad}    

\titleformat{name=\section,numberless}
  {\normalfont\bfseries}
  {}
  {0pt}
  {}

\newlength{\tabwidth}
\newlength{\tabheight}
\newlength{\tabrule}
\newlength{\tabwidthx}
\newlength{\tabheightx}

\def\gentabbox#1#2#3#4{\vbox to \tabheight{\setlength{\tabrule}{#3}%
  \setlength{\tabwidthx}{#1\tabwidth}\addtolength{\tabwidthx}{\tabrule}%

\setlength{\tabheightx}{#2\tabheight}\addtolength{\tabheightx}{-\tabheight}%
  \hbox to #1\tabwidth{%
 \hspace{-0.5\tabrule}\rule{\tabrule}{#2\tabheight}\hspace{-\tabrule}%
    \vbox to #2\tabheight{\hsize=\tabwidthx%
      \vspace{-0.5\tabrule}\hrule width\tabwidthx height\tabrule%
      \vspace{-0.5\tabrule}\vfil%
      \hbox to \tabwidthx{\hss#4\hss}%
        \vfil\vspace{-0.5\tabrule}%
      \hrule width\tabwidthx height\tabrule\vspace{-0.5\tabrule}}%
 \hspace{-\tabrule}\rule{\tabrule}{#2\tabheight}\hspace{-0.5\tabrule}}%
  \vspace{-\tabheightx}}}
\def\genblankbox#1#2{\vbox to \tabheight{\vfil\hbox to
#1\tabwidth{\hfil}}}
\def\tabbox#1#2#3{\gentabbox{#1}{#2}{0.4pt}{\strut #3}}

\catcode`\:=13 \catcode`\.=13 \catcode`\;=13
\catcode`\>=13 \catcode`\^=13
\def:#1\\{\hbox{$#1$}}
\def.#1{\tabbox{1}{1}{$#1$}}
\def>#1{\tabbox{2}{1}{$#1$}}
\def^#1{\tabbox{1}{2}{$#1$}}
\def;{\genblankbox{1}{1}\relax}
\catcode`\:=12 \catcode`\.=12 \catcode`\;=12
\catcode`\>=12 \catcode`\^=7

\newcommand\T{\mathbf{T}}

\begin{document}

\noindent
{\bf \large Orbital varieties in type $D$}

\vspace{.25in}
\noindent
WILLIAM~M.~MCGOVERN \\
\vspace{.2in}
{\it \small Department of Mathematics, Box 354350, University of Washington, Seattle, WA, 98195}
 
\vspace{.3in}
\noindent We extend the results of \cite{M20} to type $D$, again correcting the proofs of these results in \cite{M99}.  We first note that the hyperoctahedral Weyl group $W'$ of type $B_n$ acts by automorphisms on the semisimple Lie algebra $\mathfrak g$ of type $D_n$, so that the orbital variety $V(w)$ defined in \cite{M99} makes sense for any $w\in W'$; we say that $w_1,w_2\in W'$ lie in the same geometric cell if $V(w_1) = V(w_2)$ .  The classification theorem then applies to any $w\in W'$; recall that any such $w$ has well-defined left and right domino tableaux $T_L(w),T_R(w)$ of the same shape.

\newtheorem*{thm}{Theorem}
\begin{thm}
Orbital varieties contained in the orbit $\mathcal O_{\mathbf p}$ with partition $\mathbf p$ are in bijection to standard domino tableaux of shape $\mathbf p$, except that if this shape corresponds to a very even partition, then there are two orbits $\mathcal O_{\mathbf p,I},\mathcal O_{\mathbf p,II}$ attached to $\mathbf p$; standard domino tableaux of shape $\mathbf p$ such that the number $n_v$ of their vertical dominos is congruent to 0 mod 4 correspond to varieties lying in one of these orbits, say the one with Roman numeral I, while such varieties lying in the other orbit are parametrized by domino tableaux with $n_v$ congruent to 2 mod 4.  For any $w\in W'$ the tableau parametrizing the corresponding variety $V(w)$ may be obtained from the left tableau $T_L(w)$ attached by Garfinkle to $w$ by moving through certain open cycles of type $D$ or $C$ whose holes and corners lie in rows of even length, lowering the partition corresponding to the shape of the resulting tableau in the standard dominance order throughout, until this shape becomes a $D$-partition $\mathbf q$.  The orbit $\mathcal O_w$ then coincides with $\mathcal O_{\mathbf q}$ (possibly with a Roman numeral attached).  
\end{thm}

\begin{proof}
We follow the argument of \cite{M20} closely.   The operator $T_{\alpha\beta}^R$ attached to any pair $\alpha,\beta$ of adjacent simple roots with the same length carries over to this setting with the same definition.  It preserves tableau shapes except in the case $\{\alpha,\beta\} = \{e_3-e_2,e_1+e_2\}$; in that case its action is described in \cite[4.3.5]{MP16}.  We replace the operator $V_{\alpha\beta}^R$ of \cite{M20} by the operator $V_C^R$, a single-valued truncation of the operators $T_{C\alpha_1'}$ and $T_{\alpha_1'C}$ of \cite[4.4.3,4.4.10]{MP16}.  More precisely, the domain of this operator consists of all tableau pairs $(\T_1,\T_2)$ in the union of the domains of $T_{C\alpha_1'}$ and $T_{\alpha_1'C}$ such that the first four dominos of $\T_2$ form a subtableau of shape $B_1^1,B_2^1,B_3^1,B_4^1, C^1,C^2$, or $C^4$ in the notation of \cite[4.4.9]{MP16}; the operator is defined in the same way as $T_{\alpha_1'C}$ or $T_{C\alpha_1'}$, except that a tableau pair $(\T_1',\T_2')$ in its range should have either the common shape of the $\T_i$ or shape $(3,.3,1,1)$.  One further operator $V_{C'}^R$ is needed, analogous to $V_D^R$ in \cite{M20}:  it is defined on pairs $(\T_1,\T_2)$ such that the first six dominos of $\T_2$ form a subtableau of shape $(5,4,2,1)$ such that the 3- and 5-dominos are vertical.  It acts on $(\T_1,\T_2)$ by moving both $\T_i$ through the extended open cycle {\sl in type $C$} of the 6-domino and then transposing the positions of the 5- and 6-dominos in $\T_2$ in the $2\times 2$ box that they now occupy.  Calculations in types $D_4$ and $D_6$ show that all of the above operators preserve geometric left cells.  Finally we need the operator $H$, defined on pairs $(\T_1,\T_2)$ where the common shape of $\T_1$ and $\T_2$ is a single hook; it moves both tableaux through the unique open cycle of this hook and is easily seen to preserve left geometric cells.

As in \cite{M20}, we now argue as in \cite[p. 2987]{M99}, with a small modification.  Given $w\in W'$ whose left tableau $T_L(w)$ has shape that is not already a $D$-partition, if the shape of $T_L(w)$ is a single hook, then apply the operator $H$ to produce a new tableau pair whose common shape is a $D$-partition (and again a single hook).   Otherwise, we look at the largest even part $p$ occurring with odd multiplicity in the partition corresponding to this shape.  If the lowest row of length $p$ in this shape is an odd-numbered row, then it is easy to construct a tableau $\T$ having the same shape as $T_R(w)$ with the $m(m+1)$- and $(m+1)^2$-dominos having the same orientation as in $T_R(w)$, for all $m\ge3$ such that dominos with these labels appear in $T_R(w)$, and such that either the first three dominos in $\T$ form a subtableau of shape $(4,2)$ or the first four dominos of $\T$ form a subtableau of shape $(4,3,1),(4,2,2),(4,2,1,1)$, or $(3,3,2)$, as specified in the definition of $V_C^R$; the proof of the main result of \cite{M20'} shows that we may apply a sequence $\Sigma$ of operators $T_{\alpha\beta}^R$ and $V_{\alpha_1',C}^R$ to $T_R(w)$ in such a way that tableau shapes are preserved throughout and the resulting tableau is $\T$.  Applying $\Sigma$ to $(T_L(w),T_R(w)$ and then applying either the operator $T_{e_3-e_2.e_2+e_1}^R$ or $V_C^R$ , we lower the shape of both tableaux in the dominance order, arriving at $w'\in W'$ in the same geometric left cell as $w$ whose left tableau $T_L(w')$ has lower shape than that of $T_L(w)$.  If instead the lowest row of length $p$ in the shape of $T_R(w)$ is an even-numbered row, then we argue similarly, producing a tableau $\T$ in the domain of $V_C^R$ whose first five dominos form a subtableau of shape $(5,4,2,1)$ and a sequence $\Sigma$ of operators $T_{\alpha\beta}^R$ and $V_{\alpha_1',C}^R$ taking $T_R(w)$ to $\T$.  Applying $\Sigma$ to $(T_L(w),T_R(w))$ and then $V_{C'}^R$ we again get an element $w'$ in the same geometric left cell as $w$ whose tableaux $T_L(w'),T_R(w')$ have lower shape than that
of $T_L(w)$.  Iterating this process we eventually replace $w$ by another element $v$ in its geometric left cell such that the shape of $T_L(v)$ is a $D$-partition; as in \cite{M20} we see that subsequent applications of operators as above cannot lower the shape of $T_L(v)$ any further.

We now define generalized $\tau$-invariants of orbital varieties in type $D$ as was done in types $B$ and $C$ in \cite{M20}, replacing the operator $T_{\alpha\beta}^V$ occurring there for adjacent simple roots $\alpha,\beta$ of different lengths by the analogously defined operator $V_{\alpha_1',C}^V$.  Then orbital varieties in type $D$ are completely determined by their generalized $\tau$-invariants.  We then complete the classification as in \cite{M20}:  given $w,w'\in W'$ the same left tableau $T_L(w) = T_L(w')$ having the shape of a $D$-partition, any operator $T_{\alpha\beta}^V$ or $V_{\alpha_1',C}^V$ defined on $V(w)$ is also defined on $V(w')$ and preserves the shapes of $T_L(w)$ and $T_L(w')$, whence it in fact fixes these tableaux, so that $V(w)$ and $V(w')$ have the same generalized $\tau$-invariant and so must coincide.  Hence there are at most as many distinct orbital varieties as there are standard domino tableaux with shape that of a $D$-partition; but the number of such tableaux is the sum of the dimensions of the Springer representations attached to the corresponding nilpotent orbits in the ambient Lie algebra, which is in turn known to be the number of orbital varieties, so we get a natural bijection between orbital varieties and standard domino tableaux with shape that of a $D$-partition.  By \cite[p. 2986]{M99}, we know that given a $D$-partition $\mathbf p$ (possibly together with a Roman numeral) at least one $w\in W'$ has left tableau $T_L(w)$ of shape with partition $\mathbf p$ and orbital variety $V(w)$ lying in the orbit $\mathcal O_{\mathbf p}$, the bijection between standard domino tableaux and orbital varieties must send tableaux with shape $\mathbf p$ to varieties lying in $\mathcal O_{\mathbf p}$. as claimed
\end{proof}

Finally, the proof of Theorem 2 in \cite{M99} applies to show that one component of the variety $\mathcal V(w)$ of the simple highest weight module $L(w\gamma)$ of highest weight $w\gamma-\rho$, where $\gamma$ is an antidominant regular integral weight and $\rho$ is half the sum of the positive roots, is the one corresponding to the right tableau $T_R(w)$, the unique tableau of special shape obtained from the right tableau $T_R(w)$ by moving through open $D$-cycles and allowing both up and down cycles  in the sense of \cite{G93}.

\end{document}